\documentclass{article}

\def\sm{\setminus}
\def\ra{\rightarrow}
\def\sm{\setminus}
\def\ss{\subseteq}

\def\e{\epsilon}

\def\Re{\hbox{\rm Re}\,}
\def\Im{\hbox{\rm Im}\,}

\def\dbar{\overline{\partial}}

 \def\HollowBox #1#2{{\dimen0=#1 \advance\dimen0 by -#2
       \dimen1=#1 \advance\dimen1 by #2
        \vrule height #1 depth #2 width #2
        \vrule height 0pt depth #2 width #1
        \llap{\vrule height #1 depth -\dimen0 width \dimen1}%
       \hskip -#2
       \vrule height #1 depth #2 width #2}}
 \def\BoxOpTwo{\mathord{\HollowBox{6pt}{.4pt}}\;}

\def\endpf{\hfill $\BoxOpTwo$}

\def\S{{\bf S}}
\def\A{{\bf A}}

\def\bomega{\partial \Omega}
\def\sss{\subset \, \, \subset}

\font\teneufm=eufm10
\font\seveneufm=eufm7
\font\fiveeufm=eufm5
\newfam\eufmfam
\textfont\eufmfam=\teneufm
\scriptfont\eufmfam=\seveneufm
\scriptscriptfont\eufmfam=\fiveeufm

\newfam\msbfam
\font\tenmsb=msbm10  \textfont\msbfam=\tenmsb
\font\sevenmsb=msbm7  \scriptfont\msbfam=\sevenmsb
\font\fivemsb=msbm5    \scriptscriptfont\msbfam=\fivemsb
\def\Bbb{\fam\msbfam \tenmsb}

\def\RR{{\Bbb R}}
\def\CC{{\Bbb C}}

\newfam\msbbfam
\font\tenmsbb=msbm10  scaled \magstep2 \textfont\msbbfam=\tenmsbb
\font\sevenmsbb=msbm7 scaled \magstep2 \scriptfont\msbbfam=\sevenmsbb
\font\fivemsbb=msbm5  scaled \magstep2 \scriptscriptfont\msbbfam=\fivemsbb

\newtheorem{theorem}{Theorem}[section]
\newtheorem{corollary}[theorem]{Corollary}
\newtheorem{proposition}[theorem]{Proposition}
\newtheorem{lemma}[theorem]{Lemma}
\newtheorem{remark}[theorem]{Remark}

\newtheorem{definition}[theorem]{Definition}
\newtheorem{example}[theorem]{EXAMPLE}

\makeindex

\usepackage{graphicx}

\usepackage{amsmath}

\begin{document}

\begin{center}
\huge \bf
Canonical Kernels Versus Constructible Kernels\footnote{{\bf Key Words:}  canonical kernel, constructible kernel, Cauchy kernel,
Bergman kernel, Szeg\H{o} kernel}\footnote{{\bf MR Classification
Numbers:}  30C40, 32A25, 32A26, 46E22, 47B32}
\end{center}
\vspace*{.12in}

\begin{center}
\large Steven G. Krantz
\end{center}
\vspace*{.15in}

\begin{center}
\today
\end{center}
\vspace*{.2in}

\begin{quotation}
{\bf Abstract:} \sl
We study both canonical reproducing kernels and constructive reproducing
kernels for holomorphic functions in $\CC^1$ and $\CC^n$.  We compare and contrast the two, and also
develop important relations between the two types of kernels.   We prove
a new result about the relationship between these two kernels on certain
domains of finite type.
\end{quotation}
\vspace*{.25in}

\setcounter{section}{-1}

\section{Introduction}

Every working mathematician has encountered integral
reproducing formulas for holomorphic functions. The Cauchy
integral formula and the Poisson integral formula are perhaps
the two most central and important examples. These are
examples of {\it constructive} reproducing formulas (kernels)
because the integral formulas (kernels) can often be written down
explicitly or perhaps asymptotically (see [KRA2] and [KRA4]). What is of interest for
our purpose here is that there are other integral reproducing
formulas, which are canonical in nature (to be explained
below), but for which the formulas (kernels) generally {\it
cannot} be written down explicitly. Often the canonical
kernels have many attractive features, but the fact that they
are not explicit means that we do not necessarily understand
their singularities, and therefore it is difficult to analyze
them or to make estimates on them. We find it nearly
impossible to determine their mapping properties.

But there are techniques for making peace between the canonical and the constructive.
And these methods can be extremely useful in practice.  The purpose of this
paper is to explain these connections, and to explore where they might lead.
Furthermore, we shall prove a version of the Kerzman/Stein theorem
relating canonical and constructible kernels on strongly pseudoconvex domains
in the more general context of convex, finite type domains.

\section{First Principles}

In what follows, a {\it domain} is a connected open set.

It is arguable that many, if not most, constructible integral reproducing formulas are
consequences of (and often equivalent to) the fundamental theorem of calculus.
Or, in several variables, we would substitute ``Stokes's theorem'' for ``fundamental
theorem of calculus.''  As a simple example,
$$
f(x) - f(a) = \int_a^x f'(t) \, dt
$$
is an integral reproducing formula.  If $\varphi$ is a $C^1$ function
on $[a,x]$ such that $\varphi(a) = 0$ and $\varphi(x) = 1$ then we have
(from Leibniz's product rule)
$$
f(x) = \int_a^x f(t) \varphi'(t) \, dt + \int_a^x f'(t) \varphi(t) \, dt \, .
$$
In this way we have introduced the kernels $\varphi$ and $\varphi'$.
The monograph [BIN, pp.\ 103--113] contains incisive generalizations of
this last reproducing formula to the several-variable context.  It is again
worth noting that all these formulas are based on the fundamental theorem
of calculus (or Stokes's theorem).  

A perhaps more profound illustration of the connection between Stokes's theorem
and an important and central integral reproducing formula is the generalized
Cauchy integral formula.  We begin with a key but ancillary idea.
								
\begin{lemma} \sl
Let $\Omega$ be a bounded domain in $\CC$ with $C^1$ boundary.  Let
$\omega = \alpha(z) dz$ be a 1-form with coefficient that
is continuously differentiable on $\overline{\Omega}$.  Then
$$
\int_{\partial \Omega} \omega = \mathop{\int \!\!\!\int}_\Omega
  \frac{\partial \alpha}{\partial \overline{z}} d\overline{z} \wedge dz \, .
  \eqno (1.1.1)
$$
\end{lemma}

For a general 1-form $\lambda = a(z) dz + b(z) d\overline{z}$, it is useful
to write
\begin{eqnarray*}
d\lambda & = & \partial \lambda + \overline{\partial} \lambda \\
         & = & \left [ \frac{\partial a}{\partial z} \, dz \wedge dz +
	         \frac{\partial b}{\partial z} \, dz \wedge d\overline{z} \right ]
         + \left [ \frac{\partial a}{\partial \overline{z}} \, d\overline{z} \wedge dz +
	         \frac{\partial b}{\partial \overline{z}} \, d\overline{z} \wedge d\overline{z} \right ] \\
	 & = & \frac{\partial a}{\partial \overline{z}} \, d\overline{z} \wedge dz +
	           \frac{\partial b}{\partial z} \, dz \wedge d\overline{z} \, .
\end{eqnarray*}
Recall that, for a $C^1$ function $f$,
$$
\frac{\partial}{\partial z} f = \frac{1}{2} \left ( \frac{\partial}{\partial x} - i \frac{\partial}{\partial y} \right ) \, f
\qquad \hbox{and} \qquad
\frac{\partial}{\partial \overline{z}} f = \frac{1}{2} \left ( \frac{\partial}{\partial x} + i \frac{\partial}{\partial y} \right ) \, f \, .
$$
It is convenient to also introduce here the notation
$$
\partial f = \frac{\partial f}{\partial z} \ dz \qquad \hbox{and} \qquad
\overline{\partial} f = \frac{\partial f}{\partial \overline{z}} \ d\overline{z} \, .
$$

There is complex exterior differential notation for forms, and we refer the reader to [KRA2]
or [WEL] for details of that topic.

With this notation, formula (1.1.1) in the theorem becomes
$$
\int_{\partial \Omega} \omega = \mathop{\int \!\!\!\int}_\Omega \overline{\partial} \omega \, .
$$

Now we prove a generalized version of the Cauchy integral formula.
Note that it is valid for essentially {\it all} functions---not
just\index{Cauchy integral formula} holomorphic functions.

\begin{theorem}  \sl
If\ $\Omega \ss \CC$\ is a bounded domain with $C^1$ boundary and if
$f \in C^1(\overline{\Omega})$ then, for any $z \in \Omega,$
$$
f(z) = \frac{1}{2\pi i} \int_{\partial \Omega} \frac{f(\zeta)}{\zeta - z} d\zeta
                         - \frac{1}{2\pi i} \mathop{\int \!\!\!\int}_\Omega \frac{\left ( \partial f(\zeta)/\partial
\overline{\zeta} \right )}{\zeta - z}
                          d\overline{\zeta} \wedge d\zeta \, .
$$
\end{theorem}
\noindent {\bf Proof:}  \ \  \ Fix $z \in \Omega$ and choose $\epsilon > 0$ so that
$\overline{D}(z, \epsilon) \ss \Omega$.  Set $\Omega_\e = \Omega \sm \overline{D}(z,\e)$.
We apply Stokes's theorem to the form
$$
\omega(\zeta) = \frac{f(\zeta) d\zeta}{\zeta - z}
$$
on the domain $\Omega_\e$.  Note here that $\omega$ has no singularity
on $\overline{\Omega_\e}$.

Thus
$$
\int_{\partial \Omega_\e} \omega(\zeta) =
    \mathop{\int \!\!\! \int}_{\Omega_\e} \overline{\partial} \omega \, .
$$
Writing this out gives
$$
\int_{\partial \Omega} \frac{f(\zeta)}{\zeta - z} d\zeta -
   \int_{\partial D(z,\e)} \frac{f(\zeta)}{\zeta - z} d\zeta =
  \mathop{\int \!\!\! \int}_{\Omega_\e} \frac{\partial f(\zeta)/\partial \overline{\zeta}}{\zeta - z} \, d\overline{\zeta}
     \wedge d\zeta \, .   \eqno (1.2.1)
$$
Observe that we have reversed the orientation in the second integral
on the left because the disk is {\it inside} the region $\Omega$.

Now, as $\e \ra 0^+$, the integral on the right tends to
$$
\mathop{\int \!\!\! \int}_\Omega \frac{\partial f(\zeta)/\partial \overline{\zeta}}{\zeta - z} \, d\overline{\zeta} \wedge d\zeta \, .
$$
[We use here the fact that $1/(\zeta - z)$ is integrable.]
The first integral on the left does not depend on $\e$.  The second
integral on the left requires a little analysis:
$$
\int_{\partial D(z,\e)} \frac{f(\zeta)}{\zeta - z} d\zeta =
\int_0^{2\pi} \frac{f(z + \e e^{it})}{\e e^{it}} \cdot i \e e^{it} \, dt
  = i \int_0^{2\pi} f(z + \e e^{it}) \, dt \, .
$$
Now the last expression tends, as $\e \ra 0^+$, to
$2\pi i f(z)$.

Putting all this information into equation (1.2.1) yields
$$
f(z) = \frac{1}{2\pi i} \int_{\partial \Omega} \frac{f(\zeta)}{\zeta - z} \, d\zeta
    - \frac{1}{2\pi i} \mathop{\int \!\!\! \int}_{\Omega}
   \frac{\partial f/\partial \overline{\zeta}}{\zeta - z} \, d\overline{\zeta} \wedge d\zeta
$$
as desired.
\endpf
\smallskip \\

\begin{corollary}    \sl
With hypotheses as in Theorem 1.2, and the additional assumption
that $\overline{\partial} f = 0$ on $\Omega,$
we have
$$
f(z) = \frac{1}{2\pi i} \int_{\partial \Omega} \frac{f(\zeta)}{\zeta - z} d\zeta \, .
$$
\end{corollary}

\begin{remark} \rm Corollary 1.3
is the familiar Cauchy integral
formula from the analysis of one complex variable.
\end{remark}

Thus we have seen Stokes's theorem, in its complex formulation, used in
an important fashion to derive a decisive integral formula.

\section{Canonical Integral Formulas}

One of the first canonical integral formulas ever created was that of Stefan Bergman
[BER].	We shall present Bergman's idea in the context of a more general construction
due to Nachman Aronszajn [ARO].  This is the idea of a Hilbert space with
reproducing kernel.  Fortunately Aronszajn's idea also entails the Szeg\"{o} kernel
and several other important reproducing kernels.

\begin{definition} \rm
Let $X$ be any set and let ${\cal H}$ be a Hilbert space of complex-valued functions
on $X$.  We say that ${\cal H}$ is a {\it Hilbert space with reproducing kernel}
if, for each $x \in X$, the linear (point evaluation) map of the form
\begin{eqnarray*}
L_x: {\cal H} & \longrightarrow & \CC \\
            f & \longmapsto & f(x) \, ,
\end{eqnarray*}
is continuous.	 We write this as
$$
|f(x)| \leq C \cdot \|f\|_{\cal H} \, .	  \eqno (2.1.1)
$$
In this circumstance, the classical Riesz representation theorem (see [RUD])
tells us that, for each $x \in X$, there is a unique element $k_x \in {\cal H}$ such that
$$
f(x) = \langle f, k_x \rangle \qquad \forall f \in {\cal H} \, . \eqno (2.1.2)
$$

We then define a function
$$
K: X \times X \rightarrow \CC
$$
by the formula
$$
K(x,y) \equiv \overline{k_x(y)} \, .
$$
The function $K$ is the {\it reproducing kernel} for the Hilbert space ${\cal H}$.
\end{definition}

We see that $K$ is uniquely determined by ${\cal H}$ because, again by the Riesz representation
theorem, the element $k_x$ for each $x \in {\cal H}$ is unique.

It is a classical fact, and we shall not provide the details here (but see [KRA2] for all
the particulars), that the kernel $K$ may be (at least in principle) constructed by
way of a complete orthonormal basis for ${\cal H}$.  To wit, let $\{\varphi_j\}$ be
such a basis.  Then
$$
K(x,y) = \sum_{j=1}^\infty \varphi_j(x) \overline{\varphi_j(y)} \, .
$$
Here the convergence is in the Hilbert space topology in each variables. And in
fact the fundamental property $(2.1.1)$ of a Hilbert space with reproducing kernel shows
that the convergence is uniform on compact subsets of $X \times X$.

In what follows, on $\CC^n$, we use $dV$ to denote volume measure.

\begin{example} \rm
Let $\Omega$ be a bounded domain in $\CC$ or $\CC^n$.  Define
$$
A^2(\Omega) = \left \{ f \ \hbox{holomorphic on} \ \Omega: \int_\Omega |f(z)|^2 \, dV(z) < \infty \right \} \, .
$$
Clearly $A^2(\Omega)$ is a complex linear space.  It is commonly called the {\it Bergman space}.
We equip $A^2(\Omega)$ with the norm $\|f\|_{A^2(\Omega)} \equiv [\int_\Omega |f|^2 \, dV]^{1/2}$.
\end{example}

We have the following important preliminary result.

\begin{lemma} \sl
There is a constant $C = C(K, \Omega)$, depending only on
the domain $\Omega$ and on $K$ compact in $\Omega$, such that, if $f \in A^2(\Omega)$, then
$$
\sup_{z \in K} |f(z)| \leq C \left ( \int_\Omega |f(\zeta)|^2 \, dV(\zeta) \right )^{1/2} \equiv C \cdot \|f\|_{A^2(\Omega)} \, .  \eqno (2.3.1)
$$		
\end{lemma}
{\bf Proof:}  Choose $r > 0$ so that, if $z \in K$, then $B(z,r) \ss \Omega$.  Then, for $z \in K$,
we have
\begin{eqnarray*}
|f(z)| & = & \frac{1}{V(B(z,r))} \left | \mathop{\int}_{B(z,r)} f(\zeta) \, dV(\zeta) \right | \\
      & \leq & \frac{1}{\sqrt{V(B(z,r))}} \mathop{\int}_{B(z,r)} |f(\zeta)|^2 \, dV(\zeta)^{1/2}  \\
      & \leq & C(\Omega, K) \left ( \mathop{\int}_\Omega |f(\zeta)|^2 \, dV(\zeta) \right )^{1/2}  \\
      & \equiv & C(\Omega, K) \cdot \|f\|_{A^2(\Omega)} \, .		     							
\end{eqnarray*}
That proves the result.
\endpf
\smallskip  \\

With the indicated norm, the space $A^2(\Omega$ is a Hilbert space.  For the completeness, note that if
$\{f_j\}$ is a Cauchy sequence, then it will converge in the $L^2$ topology to
some limit function $g$.  But the lemma tells us that, for holomorphic functions, $L^2$ convergence
implies uniform convergence on compact sets (sometimes called {\it normal convergence}).  Hence
the limit function is holomorphic and $L^2$, thus a member of $A^2$.

Of course the lemma tells us immediately that $A^2(\Omega)$ is a Hilbert
space with reproducing kernel (just take the compact set $K$ to
be a singleton $\{z\}$).  The kernel $K$ is known as the {\it Bergman kernel}.
It is one of the most important invariants of modern function theory.

\begin{example}  \rm
Using the remark preceding Example 2.2, it is possible to calculate
the Bergman kernel for the disc.  We do so now.

Let $\psi(\zeta) = \zeta^j$ on the disc, $j = 0, 1, 2, \dots$.  The theory
of power series tells us that these functions form a complete basis for
$A^2(D)$.  Moreover, by parity, these functions are pairwise orthogonal.
A simple calculation with polar coordinates shows that
$$
\|\psi_j\|_{A^2(D)} = \sqrt{\pi/(j+1)} \, .
$$
Therefore
$$
\varphi_j(\zeta) = \sqrt{\frac{j+1}{\pi}} \zeta^j
$$
forms a complete orthonormal basis for $A^2(\Omega)$.

Using the formula preceding Example 2.2, we then see that
$$
K(z, \zeta) = \sum_{j=0}^\infty \frac{j+1}{\pi} (z \overline{\zeta})^j \, .
$$
Let $\alpha = z\overline{\zeta}$.  So we have
$$
\sum_{j=0}^\infty \frac{j+1}{\pi} \alpha^j = \frac{1}{\pi} \frac{d}{d\alpha} \sum_{j=0}^\infty a^{j+1} =
   \frac{1}{\pi} \frac{d}{d\alpha} \alpha \sum_{j=0}^\infty a^j = \frac{1}{\pi} \frac{d}{d\alpha} \alpha \cdot \frac{1}{1 - \alpha}
     = \frac{1}{\pi} \cdot \frac{1}{1 - \alpha^2} \, .
$$
In conclusion,
$$
K(z, \zeta) = \frac{1}{\pi} \cdot \frac{1}{(1 - z \overline{\zeta})^2}\, .
$$
\end{example}

\begin{example} \rm
We now want to describe the Szeg\H{o} theory.  In order to make this
work, we need to present some preliminary results about the
Bochner-Martinelli kernel and integral representation formula.

\begin{definition}   \rm
On $\CC^n$ we let
\begin{eqnarray*}
\omega(z) & \equiv & dz_1 \wedge dz_2 \wedge \dots \wedge dz_n \\
\eta(z) & \equiv & \sum_{j=1}^n (-1)^{j+1} z_j dz_1 \wedge \dots \wedge dz_{j-1} \wedge
                    dz_{j+1} \wedge \dots \wedge dz_n \, .
\end{eqnarray*}
\end{definition}
The form $\eta$ is sometimes called the {\em Leray form.}
We shall often write $\omega(\overline{z})$ to mean
$d\overline{z}_1 \wedge \dots \wedge d\overline{z}_n$ and likewise
$\eta(\overline{z})$ to mean
$\sum_{j=1}^n (-1)^{j+1} \overline{z}_j d\overline{z}_1 \wedge \dots \wedge d\overline{z}_{j-1} \wedge
                    d\overline{z}_{j+1} \wedge \dots \wedge d\overline{z}_n.$

The genesis of the Leray form
is explained by the following lemma.

\begin{lemma}    \sl
For any $z_0 \in \CC^n,$ any $\epsilon > 0,$ we have
$$
\int_{\partial B(z^0,\epsilon)} \eta(\overline{z}) \wedge \omega(z) =
         n \int_{B(z^0,\epsilon)} \omega(\overline{z}) \wedge \omega(z) \, .
$$
\end{lemma}
\noindent {\bf Proof:}  Notice that $d\eta(\overline{z}) = \overline{\partial}\eta(\overline{z})= n\omega(\overline{z}).$
Therefore, by Stokes's theorem,
$$
\int_{\partial B(z^0,\epsilon)} \eta(\overline{z}) \wedge \omega(z) =
       \int_{B(z^0,\epsilon)} d[\eta(\overline{z}) \wedge \omega(z)] \, .
$$
Of course the expression in $[\ \ \ ]$ is saturated in $dz$'s so,
in the decomposition $d = \partial + \overline{\partial},$ only
the term $\overline{\partial}$ will not die.  Thus the last line equals
$$
\int_{B(z^0,\epsilon)} [\overline{\partial}(\eta(\overline{z}))] \wedge \omega(z)
     = n \int_{B(z^0,\epsilon)} \omega(\overline{z}) \wedge \omega(z) \, . \eqno \BoxOpTwo
$$

\begin{remark}  \rm
Notice that, by change of variables,
\begin{eqnarray*} \int_{B(z^0,\epsilon)} \omega(\overline{z}) \wedge \omega(z)
           & = &   \int_{B(0,\epsilon)} \omega(\overline{z}) \wedge \omega(z)  \\
           & = & \epsilon^{2n} \int_{B(0,1)} \omega(\overline{z}) \wedge \omega(z) \, .
\end{eqnarray*}
A straightforward calculation shows that
\begin{eqnarray*}
\lefteqn{\int_{B(0,1)} \omega(\overline{z}) \wedge \omega(z)} \\
& = & (-1)^{q(n)} \cdot (2i)^n \cdot \mbox{\rm (volume of the unit ball in $\CC^n \approx \RR^{2n}$)} \, ,
\end{eqnarray*}
where $q(n) = [n(n-1)]/2.$  We denote the value of this integral by
$W(n).$
\end{remark}

\begin{theorem}[Bochner-Martinelli]   \sl
Let $\Omega \ss \CC^n$ be a bounded domain with $C^1$ boundary.
Let $f \in C^1(\overline{\Omega}).$  Then, for any $z \in \Omega,$ we have
\begin{eqnarray*}
   f(z) & = &  \frac{1}{nW(n)} \int_{\partial \Omega}
            \frac{f(\zeta)\eta(\overline{\zeta} - \overline{z}) \wedge \omega(\zeta)}{|\zeta - z|^{2n}} \\
        &   & \mbox{} - \frac{1}{nW(n)} \int_{\Omega}
              \frac{\dbar f(\zeta)}{|\zeta - z|^{2n}} \wedge \eta(\overline{\zeta} - \overline{z}) \wedge \omega(\zeta) \, .
\end{eqnarray*}
\end{theorem}
\noindent {\bf Proof:}  Fix $z \in \Omega.$  We apply Stokes's theorem to the form
$$
M_z(\zeta) \equiv \frac{f(\zeta)\eta(\overline{\zeta} - \overline{z}) \wedge \omega(\zeta)}{|\zeta - z|^{2n}}
$$
on the domain $\Omega_{z,\epsilon} \equiv \Omega \sm \overline{B}(z,\epsilon),$
where $\epsilon > 0$ is chosen so small that $\overline{B}(z,\epsilon) \ss \Omega.$
Note that Stokes's theorem does not apply to forms that have a singularity;
thus we may not apply the theorem to $L_z$ on any domain that contains
the point $z$ in either its interior or its boundary.  This observation
helps to dictate the form of the domain $\Omega_{z,\epsilon}.$  As the
proof develops, we shall see that it also helps to determine the outcome
of our calculation.

Notice that
$$
\partial (\Omega_{z,\epsilon}) = \partial \Omega \cup \partial B(z,\epsilon)
$$
but that the two pieces are equipped with opposite orientations.

\noindent Thus, by Stokes,
\begin{align}
\int_{\partial \Omega} M_z(\zeta) - \int_{\partial B(z,\epsilon)} M_z(\zeta)
& =  \int_{\partial \Omega_{z,\epsilon}} M_z(\zeta)  \notag \\
& = \int_{\Omega_{z,\epsilon}} d_\zeta(M_z(\zeta)) \, .  \tag*{(2.9.1)}  \\ \notag
\end{align}

\noindent Notice that we consider $z$ to be fixed and $\zeta$ to be
the variable.  Now
\begin{align}
d_\zeta M_z(\zeta) & = \dbar_\zeta M_z(\zeta)  \notag \\
                   & = \frac{\dbar f(\zeta) \wedge \eta(\overline{\zeta} - \overline{z}) \wedge \omega(\zeta)}{|\zeta - z|^{2n}}  \notag  \\
                   & \quad + f(\zeta) \cdot \left [ \sum_{j=1}^n \frac{\partial}{\partial \overline{\zeta}_j}
                 \left ( \frac{\overline{\zeta}_j - \overline{z}_j}{|\zeta- z|^{2n}}\right ) \right ]
                 \omega(\overline{\zeta}) \wedge \omega(\zeta) \, .   \tag*{(2.9.2)}	\\   \notag
\end{align}
Observing that
$$
\frac{\partial}{\partial \overline{\zeta}_j}
                 \left ( \frac{\overline{\zeta}_j - \overline{z}_j}{|\zeta- z|^{2n}}\right ) =
             \frac{1}{|\zeta - z|^{2n}} - n \frac{|\overline{\zeta}_j - \overline{z}_j|^2}{|\zeta - z|^{2n+2}} \, ,
$$
we find that the second term on the far right of (2.9.2) dies and we have
$$
d_\zeta M_z(\zeta) = \frac{\dbar f(\zeta) \wedge \eta(\overline{\zeta} - \overline{z}) \wedge \omega(\zeta)}{|\zeta - z|^{2n}}  \, .
$$
Substituting this identity into (2.9.1) yields
$$
\int_{\partial \Omega} M_z(\zeta) - \int_{\partial B(z,\epsilon)} M_z(\zeta)
   = \int_{\Omega_{z,\epsilon}} \frac{\dbar f(\zeta) \wedge \eta(\overline{\zeta} - \overline{z}) \wedge \omega(\zeta)}{|\zeta - z|^{2n}} \, . \eqno (2.9.3)
$$
Next we remark that
\begin{align}
\int_{\partial B(z,\epsilon)} M_z(\zeta) & = f(z) \int_{\partial B(z,\epsilon)}
             \frac{\eta(\overline{\zeta} - \overline{z}) \wedge \omega(\zeta)}{|\zeta - z|^{2n}} \notag  \\
   &  \mbox{} + \int_{\partial B(z,\epsilon)}
                  \frac{\left ( f(\zeta) - f(z) \right )\eta(\overline{\zeta} - \overline{z}) \wedge \omega(\zeta)}{|\zeta - z|^{2n}} \notag  \\
   & \equiv T_1 + T_2 \, .   \tag*{(2.9.4)}   \\ \notag
\end{align}
\noindent Since
$|f(\zeta) - f(z)| \leq C|\zeta - z|$ (and since each term of
$\eta(\overline{\zeta} - \overline{z})$ has a
factor of some $\overline{\zeta}_j - \overline{z}_j$), it follows that the integrand
of $T_2$ is of size $ O(|\zeta - z|)^{-2n+2} \approx \epsilon^{-2n + 2}.$
Since the surface over which the integration is performed has area
$\approx \epsilon^{2n-1},$ it follows that $T_2 \ra 0$ as $\epsilon \ra 0^+.$

By Lemma 2.7 and the remark following, we also have
\begin{align}
T_1 & = \epsilon^{-2n} f(z) \int_{\partial B(z,\epsilon)} \eta(\overline{\zeta} - \overline{z}) \wedge \omega(\zeta) \notag \\
    & = n\epsilon^{-2n}f(z) \int_{B(0,\epsilon)} \omega(\overline{\zeta}) \wedge \omega(\zeta) \notag \\
    & = n W(n) f(z)  \,  . \tag*{(2.9.5)} \,	 \\  \notag
\end{align}

Finally, (2.9.3)--(2.9.5) yield that
$$
\left ( \int_{\partial \Omega} M_z(\zeta) \right ) - nW(n)f(z) + o(1) =
     \int_{\Omega_{z,\epsilon}} \dbar f(\zeta) \wedge
     \left [ \frac{\eta(\overline{\zeta} - \overline{z})}{|\zeta - z|^{2n}} \right ] \wedge \omega(\zeta) \, .
$$
Since
$$
\left | \frac{\eta(\overline{\zeta} - \overline{z})}{|\zeta - z|^{2n}} \right | = O(|\zeta - z|^{-2n+1}) \, ,
$$
the last integral is absolutely convergent as $\epsilon \ra 0^+$ (remember
that $\dbar f$ is bounded).  Thus we finally have
$$
f(z) = \frac{1}{nW(n)} \int_{\partial \Omega} L_z(\zeta) -
            \frac{1}{nW(n)} \int_{\Omega} \dbar f(\zeta) \wedge
            \frac{\eta(\overline{\zeta} - \overline{z})}{|\zeta - z|^{2n}}
            \wedge \omega(\zeta)  \, .
$$
This is the Bochner-Martinelli formula.
\endpf
\smallskip \\

\begin{remark} \rm
We see that the Bochner-Martinelli formula is a quintessential
example of a constructible integral formula.  The kernel is quite
explicit, and it is the same for all domains.  For the Bergman kernel,
and for other canonical kernels that we shall see below, this latter property
does not hold.
\end{remark}

We note that the classical Cauchy integral formula in one complex variable is
an immediate consequence of our new Bochner-Martinelli formula.
	
\begin{corollary}   \sl
If $\Omega \ss \CC^n$ is bounded and has $C^1$
boundary and if $f \in C^1(\overline{\Omega})$ and $\dbar f = 0$ on $\Omega,$
then
$$
f(z) = \frac{1}{nW(n)} \int_{\partial \Omega}
                     \frac{f(\zeta) \eta(\overline{\zeta} - \overline{z})}{|\zeta - z|^{2n}} \wedge \omega(\zeta) \, .
                                \eqno (2.11.1)
$$
\end{corollary}

\begin{corollary} \sl
In complex dimension 1, the last corollary says that
$$
f(z) = \frac{1}{2\pi i} \oint_{\partial \Omega} \frac{f(\zeta)}{\zeta - z} \, d\zeta \, .
$$
\end{corollary}

Corollary 2.11 is particularly interesting.  Like the classical Cauchy formula, it gives
a constructible integral reproducing formula that is the same on all domains.  Unlike the classical
Cauchy formula, its kernel is {\it not} holomorphic in the free variable $z$.  This makes
the Bochner-Martinelli formula of limited utility in {\it constructing} holomorphic
functions.

We note that Corollary 2.11 holds for broader classes of holomorphic
functions---such as the Hardy classes.  One sees this by a simple
limiting argument.  See our discussion of $H^2$ below.

Now we turn to the development of the Szeg\H{o} theory.
Let $\Omega$ be a bounded domain in $\CC$ or $\CC^n$
with $C^1$ boundary. Define $A(\Omega)$ to be those functions which are
continuous on $\overline{\Omega}$ and holomorphic on $\Omega$. 
Identify each element of $A(\Omega)$ with its boundary trace.
Define $H^2(\Omega)$ to be the closure in the $L^2(\partial \Omega)$ norm of
$A(\Omega)$.   If $z \in \Omega$ is fixed then, by inspection of
the formula in Corollary 2.11, and the Schwarz inequality,
$$
|f(z)| \leq C \cdot \|f\|_{H^2(\Omega)} \, .
$$
Thus we see that $H^2(\Omega)$ is a Hilbert space with reproducing kernel.  We
denote the kernel for this space by $S(z, \zeta)$.

Using the comment preceding Example 2.2 of Section 2, we can actually calculate the
Szeg\H{o} kernel on the disc.  We first note that $\{z^j\}_{j=0}^\infty$ forms a
basis for the Hilbert space $H^2(D)$.  This follows from the standard theory
of power series for holomorphic functions on $D$.   It is orthogonal by parity.
It is complete by the uniqueness
of the power series expansion.  With a simple calculation, we can normalize
the basis to the complete orthonormal basis $\{1/\sqrt{2\pi} z^j\}_{j=0}^\infty$.  Thus we
see that
$$
S(z, \zeta) = \sum_{j=0}^\infty \frac{1}{2\pi} z^j \overline{\zeta}^j = \frac{1}{2\pi} \cdot \frac{1}{1 - z \cdot \overline{\zeta}} \, .
$$
Now it is instructive to write out the Szeg\H{o} integral for a function in $A(D)$:
\begin{eqnarray*}
f(z) & = & \int_0^{2\pi} f(e^{i\theta}) \cdot \frac{1}{2\pi} \cdot \frac{1}{1 - z\cdot e^{-i\theta}} \, d\theta \\
     & = & \frac{1}{2\pi i} \int_0^{2\pi} \frac{f(e^{i\theta})}{e^{i\theta} - z} \cdot i e^{i\theta} \, d\theta \\
     & = & \frac{1}{2\pi i} \oint_{\partial D} \frac{f(\zeta)}{\zeta - z} \, d\zeta \, .
\end{eqnarray*}
Thus we see that the canonical Szeg\H{o} integral formula is in fact nothing other than
the constructive Cauchy integral formula.  But only on the disc!
\end{example}

We conclude this section by noting that the integral
$$
{\bf S} g(z) = \int_{\partial \Omega} S(z, \zeta) g(\zeta) \, d\zeta
$$
defines a projection from $L^2(\partial \Omega)$ to $H^2(\Omega)$.
This is because the mapping is self-adjoint, idempotent, and fixes
$H^2$.  We call this mapping the {\it Szeg\H{o} projection}.  [Note
that the Bergman projection is constructed similarly.]

\section{Constructive Integral Formulas with Holomorphic Kernel}

In one complex variable it is easy to construct integral
formulas. The Cauchy formula is quite trivial to write down.
And it is the same for any domain.
One may also write down formulas on the ball and polydisc in
$\CC^n$. After that things become complicated.  Certainly one should
mention here the classic work [HUA] in which the Bergman kernel
is calculated for each of the Cartan bounded symmetric domains.

It was not until about 1970 that people found ways to write
down integral reproducing formulas with holomorphic kernels in
several complex variables. Here we discuss the idea.  [It should
be noted that both Bungart [BUN] and Gleason [GLE] proved some time
ago---by abstract means---that reproducing kernels that are holomorphic
in the free variable {\it always} exist.  But their methods are nonconstructive,
and the proofs quite abstract.]

Fix a non-negative integer $k$ and a strongly pseudoconvex domain
$\Omega \sss \CC^n$ with $C^{k+3}$ boundary.  Let $\rho: \CC^n \ra \RR$ be
a $C^{k+3}$ defining function for $\Omega$ with the property that
it is {\it strictly plurisubharmonic} in a neighborhood
of $\bomega.$  It is a standard fact (see [KRA2]), for which we do not provide the details, that
the function (known as the {\it Levi polynomial})
$$
L: \CC^n \times \CC^n \ra \CC
$$
given by
$$
L_P(z) = L(z,P) \equiv \rho(P) + \sum_{j=1}^n \frac{\partial \rho}{\partial z_j} (P) (z_j - P_j) \\
        \mbox{ \ \ } + \frac{1}{2} \sum_{j,k=1}^n \frac{\partial^2 \rho(P)}{\partial z_j \partial z_k}
             (z_j - P_j) (z_k - P_k)
$$
satisfies the following properties:

\begin{tabbing}
(3.1) \= \normalsize For each $P \in \CC^n,$ the function $z \mapsto L(z,P)$ is holomorphic \\
        \> (indeed, it is a polynomial);  \\
(3.2) \> For each $z \in \CC^n,$ the function $P \mapsto L(z,P)$ is $C^{k+1};$  \\
(3.3) \> For each $P \in \bomega,$ there is a neighborhood $U_P$ such that if \\
        \> $z \in \overline{\Omega} \cap \{w \in U_P: L_P(w) = 0\},$ then $z = P.$
\end{tabbing}

Our goal is to remove the need to restrict to a small neighborhood of
$P \in \bomega$ (property (3.3)) while preserving properties (3.1)--(3.3).
We proceed through a sequence of lemmas.   Following Henkin, we use the
notation
\begin{eqnarray*}
  \Omega_\delta & = & \{z \in \CC^n: \rho(z) < \delta\} ; \\
   U_\delta & = & \{z \in \CC^n: |\rho(z)| < \delta\} , \quad \delta > 0 \, .
\end{eqnarray*}
Further, let us fix the following constants:
\smallskip \\
\noindent (3.4)\ \ Choose $\delta > 0$ and $\gamma > 0$ such that
$$
 \sum_{j,k = 1}^n \frac{\partial^2 \rho}{\partial z_j \partial \overline{z}_k} (P)
           w_j \overline{w}_k \geq \gamma |w|^2 , \ \ \mbox{\rm all} \ \ P \in U_\delta, \ \ \ \mbox{\rm all}  \  w \in \CC^n \, .
$$
\noindent (3.5)\ \ Shrinking $\delta$ if necessary, we may select $\kappa > 0$ so that
$$
 |\mbox{\rm grad}\, \rho(z) | \geq \kappa \quad \mbox{for all} \  z \in U_\delta \, .
$$
\noindent (3.6)\ \ With $\delta$ as above, let
$$
 K = \sum_{|\alpha| + |\beta| \leq 3} \left \| \left ( \frac{\partial}{\partial z} \right ) ^\alpha
                              \left ( \frac{\partial}{\partial \overline{z}} \right ) ^\beta \rho(z) \right \| _{L^\infty(U_\delta)} \, .
$$

\setcounter{theorem}{6}

\begin{lemma}  \sl
There is a $\lambda > 0$ such that, if $P \in \bomega$ and
$|z - P| < \lambda$, then
$$
  2 \Re L_P(z) \leq \rho(z) - \gamma |z - P|^2 /2 \, .
$$
\end{lemma}
\noindent {\bf Proof:}  Let $\lambda = \frac{1}{2} \gamma/(K + 1).$  If $|z - P| < \lambda,$ then
$$
 \rho(z) = \rho(P) + 2 \Re L_P(z) + \sum_{j,k =1}^n \frac{\partial^2 \rho}{\partial z_j \partial \overline{z}_k} (P) (z_ j - P_j)(\overline{z}_j - \overline{P}_k) + R_P(z) \, ,
$$
where $R_P$ is the remainder term for Taylor's formula.  Therefore
\begin{align}
2 \Re L_P(z) & \leq \rho(z) - \sum_{j,k =1}^n \frac{\partial^2 \rho}{\partial z_j \partial \overline{z}_k} (P) (z_ j - P_j)(\overline{z}_j - \overline{P}_k) + |R_P(z)| \notag \\
             & \leq \rho(z) - \gamma |z - P|^2 + K|z - P|^3 \notag \\
             & \leq \rho(z) - \gamma |z - P|^2/2 \, . \tag*{$\BoxOpTwo$} \, . \\  \notag
\end{align}

\begin{lemma}  \sl
Let $\epsilon = \gamma\lambda^2/20.$
If $P \in \bomega, z \in \Omega_\epsilon, \lambda /3 \leq |z - P| \leq 2\lambda/3,$
then
$$
\Re L_P(z) < 0 \, .
$$
\end{lemma}
\noindent {\bf Proof:}  With $z$ as in the
hypotheses, we have by Lemma 3.7 that
$$
 2\Re L_P(z) \leq \epsilon - \gamma \frac{(\lambda/3)^2}{2}
                  = \gamma\frac{\lambda^2}{20} - \gamma \frac{\lambda^2}{18} < 0 \, . \eqno \BoxOpTwo
$$

\noindent We may assume that $\epsilon < \lambda < \delta < 1$
(where $\delta$ is as in (3.4) and (3.5)).  Let
$\eta: \RR \ra [0,1]$ be a $C^\infty$ function that satisfies $\eta(x) = 0$
for $x \geq 2\lambda/3$ and $\eta(x) = 1$ for $x \leq \lambda/3.$

\begin{lemma}    \sl
Fix $P \in \bomega.$  The $(0,1)$ form
$$
  f_P(z) = \left \{ \begin{array}{lcl}
              - \dbar_z\left \{\eta(|z - P|)\right \}\cdot \log L_P(z) & \mbox{\rm if} & |z - P| < \lambda , z \in \Omega_\epsilon \\
                       0 & \mbox{\rm if} & |z - P| \geq \lambda , z \in \Omega_\epsilon
                      \end{array}
             \right.
$$
is well-defined (if we take the principal branch for logarithm) and has
$C^\infty$ coefficients for $z \in \Omega_\epsilon.$  If $z$ is fixed, then
$f_P(z)$ depends $C^k$ on $P.$  Finally, $\dbar_z f_P(z) = 0$
on $\Omega_\epsilon.$  [One may
note that this
construction is valid even for $P$ sufficiently near $\partial \Omega.$]
\end{lemma}
\noindent {\bf Proof:}  On $\hbox{\rm supp}\,\{\dbar \eta(|z - P|)\}$, we have
$\lambda/3 \leq |z - P| \leq 2\lambda/3,$ so if $z$ is
also in $\Omega_\epsilon,$ then Lemma 3.8 applies and
$\Re L_P(z) < 0.$  Therefore $\log L_P(z)$ makes sense.
It follows that the form has $C^\infty$ coefficients for
$|z - P| < \lambda.$
When $|z - P| > 2\lambda/3,$ we have $\dbar_z\{\eta(|z - P|)\} \equiv 0,$
so that $f_P(z)$ is smooth.  Since $\log L_P(z)$ is holomorphic
on $\hbox{\rm supp}\, \dbar_z\eta(|z - P|),$ it follow that $\dbar f_P(z) = 0$
on all of $\Omega_\epsilon.$  The fact that $f_P(z)$ depends $C^k$ on $P$
is clear since $L_P(z)$ does. \endpf  \smallskip \\

\begin{lemma}         \sl
There is a $C^\infty$ function $u_P$ on $\Omega_\epsilon$
such that $\dbar u_P = f_P.$
\end{lemma}
\noindent {\bf Proof:}  Since $\epsilon < \delta,$ we know that $\rho_\epsilon(z) \equiv \rho(z) - \epsilon$
is a defining function for $\Omega_\epsilon;$ hence $\Omega_\epsilon$ is strongly
pseudoconvex.  By H\"{o}rmander's existence theorem for the $\overline{\partial}$ equation
(see [KRA2, Ch.\ 4]), such a function $u_P$ must exist. \endpf  \smallskip \\

We now define
$$
 \Phi(z,P) = \left \{ \begin{array}{lcl}
                      \left [\exp u_P(z)\right ] \cdot L_P(z) & \mbox{\rm if} & |z - P| < \lambda/3 \\
                      \exp\left [ u_P(z) + \eta(|z - P|) \log L_P(z) \right ] & \mbox{\rm if} & \lambda/3 \leq |z - P| < \lambda \\
                      \exp(u_P(z))                    &  \mbox{\rm if} & \lambda \leq |z - P| \, .
                        \end{array}
               \right.
$$
Notice that $\Phi$ is unambiguously defined.
To study the properties of $\Phi,$ we require two technical lemmas.

\begin{lemma}   \sl
If $U \ss \CC^n$ is any open set and $K \sss U,$ then any
$u \in C^1(U)$ satisfies
$$
 \sup_K |u| \leq C \left ( \|u\|_{L^2(U)} + \|\dbar u\|_{L^\infty(U)} \right ) \, .
$$
Here the constant $C$ depends on $U$ and $K$ but not on $u.$
\end{lemma}
\noindent {\bf Proof:}  Let $V \sss U$ be a $C^1$ domain such that $K \sss V.$
Choose $\eta \in C_c^\infty(V)$ such that $\eta \equiv 1$ on $K.$  Apply
the full Bochner-Martinelli formula to the function $\eta u$ on $V.$
Then the boundary term (the first term of the Bochner-Martinelli formula)
vanishes and the desired estimate follows directly from the integrability
of the kernel in the remaining term.

The reader who knows something about partial differential
equations will note that this result also follows from the uniform ellipticity
of the $\dbar$ operator on compact subsets of $U.$
\endpf
\smallskip \\

\begin{corollary}   \sl
 Let $\Omega \sss \CC^n$ be pseudoconvex and $K \sss \Omega.$
Let $f$ be a $\dbar-$closed $(0,1)$ form with $C^1$ coefficients.
If $u = Mf$ is the H\"{o}rmander solution to $\dbar u = f$ (see [KRA2, Ch.\ 4]), then we
have
$$
 \| u \|_{L^\infty(K)} \leq C \|f\|_{L^\infty(\Omega)} \, ,
$$
where $C$ depends only on $K$ and $\Omega$ (and not on $f$ or $u$).
\end{corollary}
\noindent {\bf Proof:}  We shall not provide the details of this
argument, but instead refer the reader to [KRA2, Ch.\ 5].
\endpf
\smallskip \\

\begin{proposition}   \sl
Assume once more that $\Omega \sss \CC^n$ has $C^{k+3}$
boundary.  Then $\Phi(\cdot,P)$ is holomorphic on $\Omega_\epsilon.$
Also there is a $C > 0,$ independent of $P,$ such that for
all $z \in \Omega_{\epsilon/2}$ we have
\medskip \\
\noindent (7.1) \ \ \ if $|z - P| < \lambda/3,$ \ \ then \ \
                   $|\Phi(z,P)| \geq C | L_P(z)|;$
\medskip \\
\noindent (7.2) \ \ \ if $|z - P| \geq \lambda/3,$ \ \ then \ \
                   $|\Phi(z,P)| \geq C.$
\end{proposition}
\noindent {\bf Proof:}  If $|z - P| \geq 2\lambda/3,$ then $\Phi(z,P) = \exp u_P(z)$
and $\dbar_z \Phi(z,P) = (\exp u_P(z)) \cdot \dbar u_P(z)
= (\exp u_P(z)) \cdot f_P(z) = 0$ by construction.
If $\lambda/3 \leq |z - P| < 2\lambda/3$ then
\begin{eqnarray*}
 \dbar_z\Phi(z,P) & = & \exp \left [ u_P(z) + \eta(|z - P|) \log L_P(z) \right ] \\
                  &   & \mbox{\ \ \ } \cdot \left \{ \dbar \left [u_P(z) + \eta(|z - P|) \log L_P(z) \right ]\right \} \\
                  & = & \exp \left [u_P(z) + \eta(|z - P|) \log L_P(z) \right ] \\
                  &   & \mbox{\ \ \ } \cdot \left \{ f_P(z) + \dbar \eta (|z - P|) \cdot \log L_P(z) \right \}
\end{eqnarray*}
since $\log L_P(z)$ is holomorphic when $\lambda/3 \leq |z - P| < 2\lambda/3.$  The last line
is $0$ by definition of $f_P.$  The calculation for $|z - P| > \lambda/3$
is trivial.  Hence we find that $\Phi_P$ is holomorphic
in the $z$ variable, $z \in \Omega_{\epsilon}.$

For the estimate, notice that $f_P$ is bounded on $\Omega_\epsilon,$ uniformly
in $P,$ so that $u_P$ is bounded on $\overline{\Omega}_{\epsilon/2}$
uniformly in $P$ (by Corollary 3.12).
So there is a $C' > 0$ such that $|\exp u_P(z)| \geq C'.$  Thus
$$
 |\Phi(z,P)| = |\exp u_P(z)| \geq C' \quad \mbox{\rm if} \ \ |z - P| \geq 2\lambda/3
$$
and
$$
  |\Phi(z,P)| = |\exp u_P(z)| \cdot |L_P(z)| \geq C' |L_P(z)| \quad \mbox{\rm if} \ \ |z - P| \leq \lambda/3 \, .
$$
For $\lambda/3 \leq |z - P| \leq 2\lambda/3,$ we have by Lemmas 3.7 and
3.8 that
\begin{eqnarray*}
 \Re L_P(z) & \leq & \frac{\epsilon}{2} - \gamma \frac{(\lambda/3)^2}{2} \\
            &  =   & \frac{\gamma \lambda^2}{40} - \gamma \frac{(\lambda/3)^2}{2} \\
            &  =   & - \frac{11 \gamma \lambda^2}{360} \, .
\end{eqnarray*}
Thus
$$
 |L_P(z)| \geq \frac{11 \gamma \lambda^2}{360} \, .
$$
We conclude that, for $\lambda/3 \leq |z - P| \leq 2\lambda/3,$
$$
 |\Phi(z,P)| \geq |\exp u_P(z)| \cdot |L_P(z)| \geq C''\frac{\gamma \lambda^2}{360} \, . \eqno \BoxOpTwo
$$

Now we would like to consider the smooth dependence of $\Phi$ on $P.$
The subtlety is that our construction of $\Phi_P$ involved solving
$\dbar u_P = f_P$ so, in principle, it appears that we must check the
smooth dependence of H\"{o}rmander's solution operator on parameters.
In fact this type of smooth dependence has been checked for various
solutions of the\ $\dbar$ problem (see [GRK1]).  But, by using a little
functional analysis, we may avoid such difficult calculations.

Now fix $z \in \Omega.$  Let $\theta \in C_c^\infty(\Omega_{\epsilon})$ satisfy
$\theta(z) = 1.$  Let $s > 2n.$  Let $M_s^\theta$ be the right inverse
to $\dbar_{0,0}$ (the H\"{o}rmander solution operator) for the pseudoconvex
domain
$\Omega_\epsilon.$  Let notation be as in (3.4) through (3.6).
Let $\mu$ be in the dual space of $W^s(\Omega,\phi_1)$ [naturally this
dual space is just $W^s(\Omega,\phi_1)$ itself].  Then
$$
  \langle \mu, M_s^\theta f_P \rangle  = \langle (M_s^\theta)^* \mu, f_P \rangle \, ,
$$
which depends $C^k$ on $P$ because $f_P$ does.

\begin{proposition}   \sl
The function $\Phi(z,P)$ depends in a $C^k$ fashion on $P$ for fixed
$z \in \Omega_\epsilon.$
\end{proposition}
\noindent {\bf Proof:}  Fix $s > 2n$ and let
$\mu = e_z$ be the point evaluation functional on
$C(\Omega) \supseteq W^s(\Omega,\phi_1).$  Then, by the preceding discussion,
$$
 \langle e_z, M_s^\theta f_P \rangle = \ \left ( M_s^\theta f_P \right ) (z) = u_P(z)
$$
depends $C^k$ on $P.$   Therefore $\Phi(z,P)$ itself depends $C^k$ on $P.$
\endpf
\smallskip \\

\begin{proposition}    \sl
Let $\Omega \ss \CC^n$ be pseudoconvex.  Let $\Omega_k = \Omega \cap \{z \in \CC^n: z_1, \dots, z_k = 0\}, k = 1,\dots,n.$
Let $A_k(\Omega) = \{f\ \mbox{holomorphic on}\ \Omega: \left. f \right |_{\Omega_k} = 0\}.$
Then there are linear operators
$$
Q_i^k: A_k(\Omega) \ra \{f\ \mbox{holomorphic on}\ \Omega\} \qquad , \quad i = 1,\dots, k \, ,
$$
such that
$$
f(z) = \sum_{i=1}^k z_i \cdot (Q_i^k f)(z)
$$
for all $f \in A_k(\Omega).$
\end{proposition}
\vskip.1in

\begin{remark}  \rm
We are primarily interested in the proposition
when $k = n.$  However the proof is by induction on $k.$
\end{remark}

Key in the proof that we are about to present is the following extension result.
We cannot provide the detailed proof of this tool, but refer
the reader to Section 5.1 of [KRA2].

\begin{theorem}  \sl
Let $\Omega \ss \CC^n$ be pseudoconvex (no assumptions about
boundary smoothness, or even boundedness, need be made).  Let
$\omega = \Omega \cap \{(z_1,\dots,z_n): z_n = 0\}.$  Let
$f: \omega \ra \CC$ satisfy the property that the map

$$
(z_1,\dots,z_{n-1}) \mapsto f(z_1,\dots,z_{n-1},0)
$$
is holomorphic on $\widetilde{\omega} = \{(z_1,\dots,z_{n-1}) \in \CC^{n-1}: (z_1,\dots,z_{n-1},0) \in \omega\}.$
Then there is a holomorphic $F: \Omega \ra \CC$ such that $\left. F \right |_{\omega} = f.$
Indeed there is a linear operator
$$
{\cal E}_{\omega,\Omega}: \{ \mbox{\rm holomorphic functions on $\omega$} \} \ra
                \{ \mbox{\rm holomorphic functions on $\Omega$} \}
$$
such that $\left. ({\cal E}_{\omega,\Omega} f) \right |_{\omega} = f.$
The operator is continuous in the topology of normal
convergence.
\end{theorem}

\paragraph{Proof of the Proposition:}  If $k = 1$ and $n$ is arbitrary then the result
follows by setting $Q_1f(z) = f(z)/z_1.$   Now suppose that the result
has been proved for $k = K - 1$ and for any $n.$

Let $\widetilde{\Omega} \equiv \{z: z_K = 0\}.$  Let $f \in A_K(\Omega).$  Then
$\widetilde{f} \equiv \left. f \right |_{\widetilde{\Omega}} \in A_{K-1}(\widetilde{\Omega}).$
Therefore, by the inductive hypothesis,
$$
\widetilde{f}(z) = \sum_{i=1}^{K-1} z_i \cdot \left (\widetilde{Q}_i^{K-1} \widetilde{f} \right ) (z) \qquad, \quad z \in \widetilde{\Omega} \, ,
$$
where $\widetilde{Q}_i^{K-1}$ are the operators assumed to exist on $\widetilde{\Omega}$
for $K - 1, n - 1.$

Now we apply Theorem 3.17.  Indeed we let
$$
Q_K^K (f)(z) = \frac{  f(z) - \sum_{i=1}^{K-1} z_i \left ( {\cal E}_{\widetilde{\Omega},\Omega} \widetilde{Q}_i^{K-1} \widetilde{f} (z) \right ) }{z_K} \, .
$$
This is well defined and holomorphic on $\Omega$ since the expression
in the numerator vanishes on $\widetilde{\Omega}.$  Also, let
$$
Q_i^K f(z) = {\cal E}_{\widetilde{\Omega},\Omega} \widetilde{Q}_i^{K-1} \widetilde{f}(z) , \quad i = 1, \dots, K-1 \, .
$$
By algebra, $f(z) = \sum_{i=1}^K z_i \cdot Q_i^K f(z),$ all $z \in \Omega.$
The induction is now complete.
\endpf  \smallskip \\

\begin{corollary}   \sl
Let $\Omega \ss \CC^n$ be pseudoconvex.  Then there are continuous
linear operators
$$
T_i: \{\mbox{holomorphic functions on}\ \Omega\} \ra \{\mbox{holomorphic functions on}\ \Omega \times \Omega\}
$$
such that, for any holomorphic $f: \Omega \ra \CC$, we have
$$
f(z) - f(w) = \sum_{i=1}^n (z_i - w_i) T_i f(z,w) , \quad \mbox{all} \ \ z, w \in \Omega \, .
$$
\end{corollary}
\noindent {\bf Proof:}  Apply Proposition 3.15
to the function $F(z,w) = f(z) - f(w)$
on the domain $\Omega \times \Omega$ with coordinates
\begin{eqnarray*}
       z'_1 & = & z_1 - w_1 \\
            & \vdots &      \\
       z'_n & = & z_n - w_n \\
       z'_{n+1} & = & z_1   \\
            & \vdots &      \\
       z'_{2n}  & = & z_n \, .
\end{eqnarray*}
The continuity will follow from the closed graph theorem.
\endpf  \smallskip \\

\begin{proposition}[Hefer's Lemma]      \sl
Let $\Omega \ss \CC^n$ be strongly pseudoconvex with
$C^4$ boundary.  Let $\Phi: \Omega_{\e/2} \times \bomega \ra \CC$ be the
$C^1$ singular function constructed above.  Then we may
write
$$
\Phi(z,\zeta) = \sum_{i=1}^n (\zeta_i - z_i) \cdot P_i(z,\zeta) , \quad z \in \Omega_{\e/2} , \quad \zeta \in \bomega \, ,
$$
where each $P_i$ is holomorphic in $z \in \Omega_{\e/2}$ and $C^1$ in $\zeta \in \bomega.$
\end{proposition}
\noindent {\bf Proof:}  Fix $\zeta \in \bomega.$  Apply Corollary 3.18 to the function
$\Phi_\zeta(\cdot) = \Phi(\cdot,\zeta)$ on $\Omega_{\e/2}.$  So
$$
\Phi(z,\zeta) - \Phi(w,\zeta) = \sum_{i=1}^n (z_i - w_i) [(T_i\Phi_\zeta)(z,w)] \, .
$$
Since this is true for all $w \in \Omega_{\e/2},$ we may set $w = \zeta \in \bomega$
to obtain
$$
\Phi(z,\zeta) = \sum_{i=1}^n (z_i - \zeta_i)\left [ (T_i \Phi_\zeta)(z,\zeta) \right ] \equiv
                   \sum_{i=1}^n (\zeta_i - z_i) P_i(z,\zeta) \, .
$$
It remains to check that $P_i$ is $C^1$ in $\zeta.$  For this, it is enough
to verify that $(T_i\Phi_\zeta)(z,w)$ is $C^1$ in $\zeta.$  But, just as
in the proof of Proposition 3.14, we let $e_{(z,w)}$ be the point evaluation
functional on $\Omega \times \Omega$ and observe that
$$
(T_i \Phi_\zeta)(z,w) = \langle e_{(z,w)}, T_i \Phi_\zeta \rangle = \langle T_i^* e_{(z,w)}, \Phi_\zeta\rangle \, .
$$
The last expression is $C^1$ in $\zeta$ by Proposition 3.14.
\endpf  \smallskip \\

We now quickly review the Cauchy-Fantappi\'{e} formula.  See [KRA2, Ch.\ 1] for
the details.

\begin{theorem} \sl
Let  $\Omega \sss \CC^n$ be a domain
with $C^1$ boundary.  Let $w(z,\zeta) = (w_1(z,\zeta),\dots,w_n(z,\zeta))$
be a $C^1,$ vector-valued function on 
$\bar{\Omega} \times \bar{\Omega} \sm \{\mbox{\rm diagonal}\}$ that satisfies
$$ 
\sum_{j=1}^n w_j(z,\zeta) (\zeta_j - z_j) \equiv 1 . 
$$
Then, using the notation from Section 2, we have for any
$f \in C^1(\bar{\Omega}) \cap \{\mbox{\rm holomorphic functions on $\Omega$}\}$
and any $z \in \Omega$ the formula
$$  
 f(z) = \frac{1}{nW(n)} \int_{\partial \Omega} f(\zeta) \eta(w) \wedge \omega(\zeta) . 
$$
\end{theorem}

We see that the Cauchy-Fantappi\'{e} formula is a direct generalization
of the Bochner-Martinelli formula discussed above.  Now we can give
the punchline of this development.

\begin{theorem}[Henkin \mbox{[2]}]    \sl
Let $\Omega \ss \CC^n$ be a strongly pseudoconvex domain with $C^4$ boundary.  Let
$\Phi: \Omega_{\e/2} \times \bomega \ra \CC$ be the Henkin singular
function.  Define
$$
w_i(z,\zeta) = \frac{P_i(z,\zeta)}{\Phi(z,\zeta)} \ , \ i = 1,\dots, n \, .
$$
Here $P_i(z,\zeta)$ are as in Proposition 3.19.  Just as in our earlier discussion, let
$$
\eta(w) = \sum_{i=1}^n (-1)^{i+1} w_i dw_1 \wedge \cdots \wedge dw_{i-1} \wedge dw_{i+1} \wedge \cdots \wedge dw_n
$$
and
$$
\omega(\zeta) = d\zeta_1 \wedge \cdots \wedge d\zeta_n \, .
$$
Then, for any $f \in C^1(\overline{\Omega}) \cap \{\mbox{holomorphic functions on} \ \Omega\},$
we have the integral representation
$$
f(z) = \int_{\bomega} f(\zeta) \eta(w) \wedge \omega(\zeta) \,.
$$
\end{theorem}
\noindent {\bf Proof:}  The functions $w_i$ satisfy
$$
\sum_{i=1}^n w_i(z,\zeta) (\zeta_i - z_i) \equiv 1 , \quad z \in \Omega , \quad \zeta \in \bomega \, .
$$
Now apply the Cauchy-Fantappi\`{e} formula (see [KRA2, Ch.\ 5]). \endpf
\smallskip \\

\begin{corollary}  \sl
With notation as in the theorem, we have
$$
f(z) = \int_{\bomega} f(\zeta) \frac{K(z,\zeta)}{\Phi^n(z,\zeta)} d\sigma(\zeta) \, , \eqno (3.22.1)
$$
where $K: \Omega_{\e/2} \times \bomega$ is holomorphic in $z$ and continuous
in $\zeta.$  In fact, $K(z,\zeta) d\sigma(\zeta) = \eta(z) \wedge \omega(\zeta).$
\end{corollary}
\noindent {\bf Proof:}  See [KRA2, Ch.\ 5].
\endpf  \smallskip \\

Of course Corollary 3.22 gives us a constructive integral reproducing formula
with kernel that is holomorphic in the free $z$ variable.  This is a very
useful device, and important for the function theory of several complex variables.

\section{Asymptotic Expansion for the Canonical Kernel}

C. Fefferman [FEF] made an important contribution in 1974 when he produced an asymptotic
expansion for the Bergman kernel of a strongly pseudoconvex domain.  Basically he
was able to write
$$
K(z,\zeta) = P(z,\zeta) + {\cal E}(z, \zeta) \, ,
$$
where $P$ (the principal term) is, in suitable local coordinates, the Bergman kernel of the ball
and ${\cal E}$ (the error term) is a term of strictly lower order (in some measurable sense).
This powerful formula gives one a means for calculating mapping properties of the Bergman integral.
Fefferman himself used the formula to calculate the boundary asymptotics of Bergman metric
geodesics (for the purpose of proving the smooth boundary extension of biholomorphic mappings).  Fefferman
states in his paper---although the details have never been worked out---that there is a similar
asymptotic expansion for the Szeg\H{o} kernel of a strongly pseudoconvex domain.

At about the same time, Boutet de Monvel and Sj\"{o}strand [BOS] used the technique of Fourier
integral operators [HOR] to directly derive an asymptotic expansion for the Szeg\H{o} kernel of
a strongly pseudoconvex domain.  This expansion is quite similar to Fefferman's:  there is
a principal term, which in suitable local coordinates is the Szeg\H{o} kernel of the ball, and
there is an error term which is of lower order.  It is not known whether the techniques of [BOS] can
be used to derive an asymptotic expansion for the Bergman kernel.

The main purpose of the present paper is to consider another method, due to Kerzman and Stein, for deriving asymptotic
expansions for the canonical kernels that is more elementary and uses less machinery.
Fefferman's rather complicated argument uses Kohn's solution of the $\overline{\partial}$-Neumann
problem as well as the theory of nonisotropic singular integrals.  Boutet de Monvel and Sj\"{o}strand's
argument uses the theory of Fourier integral operators.	  The method of Kerzman and Stein [KES] that
we treat here uses only basic complex function theory and a little functional analysis.

At this time there are virtually no results about asymptotic expansions for the canonical kernels
on weakly pseudoconvex domains.  Some interesting partial results appear in [HAN].   But
see our Section 6 below.

\section{The Relation Between Constructive Kernels and Canonical Kernels on Strongly Pseudoconvex Domains}

The ideas that we present now have thus far only been developed on
strongly pseudoconvex domains. It is an important open problem
to determine how to carry out a similar program on finite type
domains or more general domains.

In previous sections, we have defined the Szeg\H{o} projection ${\bf S}: L^2(\partial \Omega) \ra H^2(\Omega)$.
We also have a mapping ${\bf H}: L^2(\partial \Omega) \ra H^2(\Omega)$ that is determined by
the Henkin kernel of Corollary 3.22.   We note that {\bf H} defines a bounded operator from $L^2(\partial \Omega)$
to $H^2(\Omega)$ (the Hardy space---see [KRA, Chapter 8]) for the following reason.  

It is known that $\partial \Omega$, when equipped with balls coming from the complex structure
and the usual boundary area measure (see [NSW1], [NSW2]), is a space of homogeneous type in the
sense of Coifman and Weiss [COW].  Further,  it is straightforward to verify that the Henkin operator {\bf H} satisfies the 
hypotheses of the David-Journ\'{e} $T1$ theorem for spaces of homogeneous type (see [CHR] for a nice
exposition of these ideas).  Thus we may conclude that the Henkin operator maps $L^2(\partial \Omega)$ to
$L^2(\partial \Omega)$.  
Since the Henkin kernel also
obviously maps $L^2(\partial \Omega)$ to holomorphic functions, we may conclude
that the Henkin integral maps $L^2(\partial \Omega)$ to $H^2(\Omega)$.

Now of course ${\bf S}$, being a projection, is self-adjoint.  So ${\bf S} = {\bf S}^*$.
It is not at all true that ${\bf H} = {\bf H}^*$, but one may calculate (see below for the
details) that ${\bf A} \equiv {\bf H}^* - {\bf H}$ is small in a measurable sense.

We also have
$$
{\bf H} \S = \S \ , \qquad \quad \S {\bf H}^* = \S \ ,
$$
$$
\S {\bf H} = {\bf H} \ , \qquad \quad {\bf H}^* \S = {\bf H}^* \, .
$$
Let us discuss these four identities for a moment.

For the first, notice that $\S$ is the projection onto $H^2$, and ${\bf H}$ preserves holomorphic
functions.  So certainly ${\bf H}\S = \S$.  For the second, we calculate that
$$
\langle \S {\bf H}^* x, y \rangle = \langle {\bf H}^* x, \S y\rangle = \langle x, {\bf H} \S y\rangle
= \langle x, \S y \rangle
$$
(because ${\bf H}$ preserves holomorphic functions) and
thus $= \langle \S x, y \rangle$.  Hence $\S {\bf H}^* = \S$.  For the third, notice that ${\bf H}$ maps
to the holomorphic functions, and $\S$ preserves holomorphic functions.  And, for the fourth,
we calculate that
$$
\langle {\bf H}^* \S x, y\rangle = \langle \S x, {\bf H} y \rangle = \langle x, \S{\bf H} y\rangle =
   \langle x, {\bf H} y\rangle = \langle {\bf H}^* x, y \rangle \, .
$$
In conclusion, ${\bf H}^* \S = {\bf H}^*$.

Now we see that
$$
\S\A = \S({\bf H}^* - {\bf H}) = \S{\bf H}^* - \S{\bf H} = \S - {\bf H} \, .
$$
As a result,
$$
\S = {\bf H} + \S \A
$$
so
$$
\S({\bf I} - \A) = {\bf H} \, .
$$
In conclusion,
$$
\S = {\bf H} ({\bf I} - \A)^{-1}  \, .
$$

If indeed we can show that $\A$ is norm small in a suitable sense, then
$({\bf I} - \A)^{-1}$ is well defined by a Neumann series.
Thus we may write
$$
\S = \bf H + {\bf H} \A + {\bf H} \A^2 + \cdots + {\bf H}\A^j + {\bf H}\A^{j+1} + \cdots \, .
$$
Hence we have expressed the Szeg\H{o} projection $\S$ as an asymptotice expansion
in terms of the Henkin projection ${\bf H}$.  By applying this asymptotic expansion
to the Dirac delta mass, this last formula
can be translated into saying that the Szeg\H{o} {\it kernel} $S$ can be
written as an asymptotic expansion in terms of the Henkin {\it kernel}.

It should be noted that Ewa Ligocka [LIG] has shown that these same ideas
may be applied to expand the Bergman kernel in an asymptotic expansion
in terms of the Henkin kernel.  We shall not treat the details of
her argument here.

\section{A version of the Kerzman/Stein Theorem on Finite Type Domains}

Here we treat a version of the main result of [KES] on certain
finite type domains.  Since most of the steps follow [KES] rather
closely, we shall outline much of the proof.

We shall concentrate in this section on a convex domain $\Omega
= \{z: \rho(z) < 0\}$ of finite type in $\CC^2$ (see [KRA2] for a
thorough discussion of the concept of finite type). In
particular, the construction of Henkin's reproducing formula
with holomorphic kernel as in Section 3 is straightforward for
such domains. For one can let $\Phi(z,P)$, for $P \in \partial
\Omega$, be given by
$$
\Phi(z,P) = \rho(P) + \sum_{j = 1}^n \frac{\partial \rho}{\partial \zeta_j} (P) (z_j - P_j) \, .
$$
No construction of $\Phi$ using the $\overline{\partial}$ problem, as in Lemmas 3.10 and 3.12 and
the discussion adhering thereto, is needed.  And we also have no need for
the quadratic terms that occur in the Levi polynomial as discussed at the beginning
of Section 3.

As a consequence, there is no need for our Lemma 3.7 in this context, and Proposition 3.14
is automatic.  Proposition 3.15 and Corollary 3,18 are true for all pseudoconvex domains,
and Proposition 3.19 is immediate from our new definition of $\Phi$.  Thus we derive
a suitable version of Theorem 3.21 and Corollary 3.22 in the present context.   We note that Range
[RAN1], [RAN2] has explored these ideas on domains of this type.

Now, looking at [KES], we see that a crucial ingredient is the function $g$ which
is defined on page 202.  In our context, the function $g$ is simply
$$
g(z,\zeta) = \rho(\zeta) + \sum_{j=1}^n \frac{\partial \rho}{\partial z_j} (z) (z_j - \zeta_j) \, .
$$
Then the key properties of $g$ enunciated in (2.1.6)--(2.1.10)of [KES]  are immediate.  Also,
for us, the kernel of the operator {\bf A} will have the same form as in the paper [KES] of
Kerzman and Stein.  The only key difference is in the estimate from below on $|g(z, \zeta)|$.
Whereas, for Kerzman and Stein, the estimate is
$$
|g(z, \zeta) | \geq C \bigl [\hbox{dist}(z, \partial \Omega) + c |\zeta - z|^2 \bigr ] \, ,    \eqno (6.1)
$$
for us it is
$$
|g(z, \zeta) | \geq C \bigl [\hbox{dist}(z, \partial \Omega) + c |\zeta - z|^{2m} \bigr ] \, .  \eqno (6.2)
$$
Here $m$ comes from the inequality
$$
\rho(\zeta + \tau) \geq C |\tau|^{2m}	\eqno (6.3)
$$
for $\tau$ a complex tangent vector at $\zeta \in \partial \Omega$.   And it is plain that $(6.3)$ holds
because the domain $\Omega$ is of finite type in $\CC^2$.  These ideas are discussed in detail
in [RAN1] and [RAN2].

It is plain to see that the estimate in $(6.2)$ is not as favorable as the one $(6.1)$ that Kerzman and Stein
got to deal with in the strongly pseudoconvex case.  But we have the advantage that the boundary
of our domain $\Omega$ has a different geometry.  For any given boundary point $\zeta \in \partial \Omega$,
we have the inequality $(6.3)$, which says that the boundary is $m^{\rm th}$-order flat at $\zeta$.
Thus one can analyze the expressions
$$
\frac{\phi(z)}{[g(\zeta, z)]^n}
$$
just as in [KES], and use the methods of [KRA1], to find that the kernel of $H$ is of weak
type $1$ in the $z$ variable, for each fixed value of $\zeta$.  And also of weak type
$1$ in the $\zeta$ variable, for each fixed value of $z$.

It follows then, because of (2.1.8) and (2.2.2) in [KES], that the kernel of {\bf A} is in fact of weak
type $1 + \delta$, some $\delta >0$, in the variable $z$ for each fixed
value of $\zeta$.  And it is of weak type $1 + \delta$ in the variable $\zeta$,
for each fixed value of $z$ (it can also be useful here to use nonisotropic polar
coordinates--see [FOS]).  Thanks to a nice lemma of Folland and Stein (which is
presented in detail in [KRA3, Theorem 9.7.7]), we then know that the operator {\bf A} maps $L^p$ to $L^{p + \epsilon}$, some
$\epsilon > 0$ (where $\epsilon$ depends in an explicit fashion on $\delta$).

But, more importantly, we can make the following analysis.
For $\lambda > 0$, we can define a fractional differential operator of order $1 - \lambda$
by
$$
{\cal D}^\lambda f (x) = (- \bigtriangleup )^{1/2} \left ( \int \frac{f(t)}{|x - t|^{N - (1-\lambda)}} \, dt \right ) \, .
$$
Here $\bigtriangleup$ is the usual Euclidean Laplacian, and the fractional power is calculated
using the Fourier transform.
Correspondingly, we let the fractional integration operator of order $\mu$ be given by
$$
{\cal I}^\mu f(x) = \int \frac{f(t)}{|x - t|^{N - \mu}} \, dt \, .
$$
Of course, in practice, we would define these operators on a coordinate patch in $\partial \Omega$.

It is easy to check, using elementary Fourier analysis, that ${\cal I}^\mu \circ {\cal D}^\lambda = \hbox{id}$ \
preciseley when $\lambda = \mu$ (and of course $0 < \lambda < 1$, $0 < \mu < 1$).

Now we can write
$$
{\bf A} = {\cal I}^\lambda \circ {\cal D}^\lambda \circ {\bf A} = {\cal I}^\lambda \circ \bigl ( {\cal D}^\lambda \circ {\bf A} \bigr ) \, .
$$
Here $\lambda$ is chosen to be quite small relative to $\epsilon$.  Denote the operator in parentheses on
the right of this last formula by ${\cal T}^\lambda$.  So
$$
{\bf A} = {\cal I}^\lambda \circ {\cal T}^\lambda \, .
$$
If $\lambda$ is small enough and positive, then the kernel of ${\cal T}^\lambda$ will be of weak
type $1 + \delta/2$ in each variable.  So it will certainly map $L^p$ {\it at least}
to $L^p$, and in fact to a higher-order Lebesgue space.  And certainly
${\cal I}^\lambda$ will map $L^2$ to a Sobolev space
of positive order.  As a result (and Kerzman/Stein discuss this point in
some detail in their paper), by Rellich's lemma, ${\bf A}$ is a compact operator on $L^2$.

It follows then that the asymptotic expansion of $S$ in terms of $H$ given
by formula $(*)$ is a valid, convergent expansion.  And that is what we wished
to prove.

We conclude by noting that the paper [KES] makes decisive use of the Heisenberg group
approximation technique introduced in [FOS] in order to obtain these last results
about the mapping properties of {\bf A}.  In the finite type case, one might consider
using the ideas in [ROS] to imitate those arguments.  We have taken a different approach
to the matter in the present paper.

\section{Concluding Remarks}

It is a matter of definite interest to be able to expand the canonical, but
non-explicit, Szeg\H{o} kernel in terms of a more explicit kernel like the
Henkin kernel.  Until now, this expansion had only been achieved on
strongly pseudoconvex domains.  We have indicated here a way to perform
the procedure on a class of finite type domains in $\CC^2$.  One might
anticipate that the result ought to be true on all finite type domains in
any dimension.  Of course the analytic difficulties attendant to obtaining
such a result are rather formidable.

We hope to explore the matter in future papers.

\newpage

\noindent {\Large \sc References}
\bigskip  \\

\begin{enumerate}						

\item[{\bf [ARO]}] N. Aronszajn, Theory of reproducing kernels, {\it Trans.
Am. Math. Soc.} 68(1950), 337-404.				     	

\item[{\bf [BER]}]  S. Bergman,  \"{U}ber die Entwicklung der harmonischen Funktionen
der Ebene und des Raumes nach Orthogonal funktionen, {\it Math.\ Annalen} 96(1922), 237--271.

\item[{\bf [BIN]}] O. V. Besov, V. P. Ilin, and S. M. Nikol'skii, {\it
Inegral Representations of Functions and Imbedding Theorems}, V. H.
Winston, New York 1978--1979.

\item[{\bf [BOS]}] L. Boutet de Monvel and J. Sj\"{o}strand, Sur la
singularit\'{e} des noyaux de Bergman et Szeg\"{o}, {\em Soc. Mat. de
France Asterisque} 34-35(1976), 123--164.
				
\item[{\bf [BUN]}]  L. Bungart, Holomorphic functions with values in locally convex
spaces and applications to integral formulas, {\it Trans. Am. Math. Soc.}
111(1964), 317-344.

\item[{\bf [CHR]}] M. Christ, {\it Lectures on Singular Integral
Operators}, CBMS Regional Conference Series in Mathematics, 77,
Published for the Conference Board of the Mathematical
Sciences, Washington, DC; by the American Mathematical
Society, Providence, RI, 1990.

\item[{\bf [COW]}] R. R. Coifman and G. Weiss, {\it Analyse
harmonique non-commutative sur certains espaces
homog\'{e}nes}, (French) \'{E}tude de certaines int\'{e}grales
singuli\'{e}res, Lecture Notes in Mathematics, Vol. 242,
Springer-Verlag, Berlin-New York, 1971.

\item[{\bf [FEF]}] C. Fefferman, The Bergman kernel and biholomorphic
mappings of pseudoconvex domains, {\em Invent. Math.} 26(1974), 1--65.

\item[{\bf [FOS]}] G. B. Folland and E. M. Stein, Estimates for
the $\dbar_b$ complex and analysis on the Heisenberg group,
{\em Comm. Pure Appl. Math.} 27(1974), 429-522.

\item[{\bf [GLE]}]  A. Gleason, The abstract theorem of Cauchy-Weil, 
{\it Pac. J. Math.} 12(1962), 511-525.

\item[{\bf [GRK1]}] R. E. Greene and S. G. Krantz, Deformation
of complex structures, estimates for the $\dbar$ equation, and
stability of the Bergman kernel, {\em Adv. Math.} 43(1982),
1--86.

\item[{\bf [HAN]}] N. Hanges, Explicit formulas for the Szeg\"{o} kernel
for some domains in $\CC^2$. {\it J. Funct.\ Anal.} 88 (1990), 153--165.

\item[{\bf [HEN]}] G. M. Henkin, Integral representations of
functions holomorphic in strictly pseudoconvex domains and
some applications, {\em Mat. Sb.} 78(120)(1969), 611-632; {\em
Math. U.S.S.R. Sb.} 7(1969), 597-616.

\item[{\bf [HOR]}] L. H\"{o}rmander, Fourier integral operators. I, {|it
Acta Math.} 127(1971), 79--183.

\item[{\bf [HUA]}] L. K. Hua, {\it Harmonic Analysis of
Functions of Several Complex Variables in the Classical
Domains}, American Mathematical Society, Providence, 1963.
				
\item[{\bf [KES]}] N. Kerzman and E. M. Stein, The Szeg{\bf H}{o}
kernel in terms of Cauchy-Fantappi\'{e} kernels, {\it Duke
Math.\ Journal} 25(1978), 197--224.

\item[{\bf [KRA1]}] S. G. Krantz, Optimal Lipschitz and $L^{p}$
regularity for the equation $\overline{\partial} u = f$ on strongly
pseudo-convex domains, {\it Math. Annalen} 219(1976), 233-260.
						
\item[{\bf [KRA2]}] S. G. Krantz, {\it Function Theory of
Several Complex Variables}, $2^{\rm nd}$ ed., American
Mathematical Society, Providence, RI, 2001.

\item[{\bf [KRA3]}] S. G. Krantz, {\it Explorations in Harmonic
Analysis, with Applications in Complex Function Theory and the
Heisenberg Group}, Birkh\"{a}user Publishing, Boston, 2009.

\item[{\bf [KRA4]}] S. G. Krantz, Calculation and estimation of
the Poisson kernel, {\it J. Math.\ Anal.\ Appl.}
302(2005)143--148.

\item[{\bf [LIG]}] E. Ligocka, The H\"{o}lder continuity of the Bergman
projection and proper holomorphic mappings, {\it Studia Math.} 80(1984),
89--107.

\item[{\bf [NSW1]}] A. Nagel, E. M. Stein, and S. Wainger,
Boundary behavior of functions holomorphic in domains of
finite type, {\it Proc. Nat. Acad. Sci. USA} 78(1981),
6596-6599.
				 
\item[{\bf [NSW2]}]  A. Nagel, E. M. Stein, and S. Wainger,
Balls and metrics defined by vector fields, I, Basic properties,
{\it Acta Math.} 155(1985), 103--147. 

\item[{\bf [RAN1]}]  R. M. Range, H\"{o}lder estimates for
$\overline{\partial}$ on convex domains in $\CC^2$ with
real analytic boundary, Several complex variables ({\it Proc.\
Sympos.\ Pure Math.}, Vol. XXX, Part 2, Williams College,
Williamstown, Mass., 1975), pp.\ 31--33. Amer.\ Math.\ Soc.,
Providence, R.I., 1977.

\item[{\bf [RAN2]}] R. M. Range, On H\"{o}lder
estimates for $\overline{\partial} u=f$ on weakly pseudoconvex
domains, Several Complex Variables (Cortona, 1976/1977), pp.\
247--267, Scuola Norm.\ Sup.\ Pisa, Pisa, 1978.

\item[{\bf [ROS]}] L. P. Rothschild and E. M. Stein, L. P. Rothschild and
E. M. Stein Hypoelliptic differential operators and nilpotent groups, {\it
Acta Math.} 137(1976), 247-320.

\item[{\bf [RUD]}] W. Rudin, {\it Real and Complex Analysis},
McGraw-Hill, New York, 1966.

\item[{\bf [WEL]}]  R. O. Wells, {\it Differential Analysis on Complex
Manifolds}, Springer-Verlag, New York, 1979.

\end{enumerate}
\vspace*{.17in}

\begin{quote}
Department of Mathematics \\
Washington University in St.\ Louis  \\
St.\ Louis, Missouri 63130 \ \ U.S.A.  \\
{\tt sk@math.wustl.edu}
\end{quote}

\end{document}